\documentclass{amsart}
\usepackage{amssymb} 
\usepackage{amsmath} 
\usepackage{amscd}
\usepackage{amsbsy}
\usepackage{comment}
\usepackage[matrix,arrow]{xy}
\usepackage{hyperref}
\usepackage{enumitem}

\DeclareMathOperator{\Norm}{Norm}

\DeclareMathOperator{\GL}{GL}

\DeclareMathOperator{\Gal}{Gal}

\DeclareMathOperator{\ord}{ord}

\DeclareMathOperator{\lcm}{lcm}

\DeclareMathOperator{\Aut}{Aut}

\newcommand{\Q}{{\mathbb Q}}
\newcommand{\Z}{{\mathbb Z}}

\newcommand{\F}{{\mathbb F}}
\newcommand{\cA}{\mathcal{A}}
\newcommand{\cE}{\mathcal{E}}
\newcommand{\cB}{\mathcal{B}}
\newcommand{\cC}{\mathcal{C}}
\newcommand{\cD}{\mathcal{D}}

\newcommand{\cI}{\mathcal{I}}

\newcommand{\ga}{{\mathfrak{a}}}

\newcommand{\ii}{\mathbf{i}}
\newcommand{\xx}{\mathbf{x}}
\newcommand{\yy}{\mathbf{y}}

\def\mod#1{{\ifmmode\text{\rm\ (mod~$#1$)}
\else\discretionary{}{}{\hbox{ }}\rm(mod~$#1$)\fi}}

\begin {document}

\newtheorem{thm}{Theorem}
\newtheorem{lem}{Lemma}[section]
\newtheorem{prop}[lem]{Proposition}

\newtheorem*{conj}{Conjecture}

\theoremstyle{definition}

\theoremstyle{remark}

\title[]{A conjecture of Erd\H{o}s, supersingular primes and short character sums}

\author[Michael Bennett]{Michael A. Bennett}
\address{Department of Mathematics, University of British Columbia, Vancouver, B.C., V6T 1Z2 Canada}
\email{bennett@math.ubc.ca}

\author{Samir Siksek}
\address{Mathematics Institute, University of Warwick, Coventry CV4 7AL, United Kingdom}
\email{S.Siksek@warwick.ac.uk}
\thanks{
The first author
was supported by NSERC, while the second author was supported by  an 
EPSRC {\em LMF: L-Functions and Modular Forms} Programme Grant
EP/K034383/1.
}

\date{\today}

\begin {abstract}
If $k$ is a sufficiently large positive integer, we show that the Diophantine equation
$$
n ( n+d) \cdots (n+ (k-1)d ) = y^{\ell}
$$
has at most finitely many solutions in positive integers $n, d, y$ and $\ell$, with $\gcd (n,d)=1$ and $\ell \geq 2$.
Our proof relies upon Frey-Hellegouarch curves and results on supersingular primes for elliptic curves without complex multiplication, derived from upper bounds for 
short character sums and sieves, analytic and combinatorial.
\end {abstract}
\maketitle

\section{Introduction} \label{intro}

In 1975, Erd\H{o}s and Selfridge \cite{ErSe} solved a long-open problem, originally posed by Liouville \cite{Lio} in 1857, 
proving that the product of two or more consecutive nonzero integers can never
be a perfect power:

\begin{thm}[Erd\H{o}s - Selfridge, 1975] \label{ES}
 The Diophantine equation
\begin{equation} \label{start}
n (n+1) \cdots (n+k-1) = y^{\ell}
\end{equation}
has no solutions in positive integers $n, k, y$ and $\ell$ with $k, \ell \geq 2$.
\end{thm}

The proof, rather surprisingly, relies upon  a combination of clever elementary
and graph theoretic arguments. Earlier work on equation (\ref{start}), from
Liouville onwards, had either depended upon results from multiplicative number
theory or upon Diophantine approximation (as, for example, in oft-cited but
unpublished work of Erd\H{o}s and Siegel where a result similar to Theorem
\ref{ES} was obtained for suitably large $n$).

An apparently rather more difficult problem is to derive an analogue of Theorem
\ref{ES} for products of consecutive terms in arithmetic progression, 
and this is the subject of the following famous conjecture, 
widely attributed to Erd\H{o}s (see for
example \cite{ShTi}):
\begin{conj}(Erd\H{o}s)
There is a constant $k_0$ such that the Diophantine equation
\begin{equation} \label{main-eq}
n(n+d)(n+2d) \cdots (n+(k-1)d)=y^\ell, \qquad \gcd(n,d)=1
\end{equation}
has no solutions in positive integers $n$, $d$, $k$, $y$, $\ell$, with $\ell \ge 2$ and $k \ge k_0$.
\end{conj}
Without the condition $\gcd(n,d)=1$ it is easy to construct a plethora
of artificial solutions. As pointed out by Erd\H{o}s
and Selfridge, equation \eqref{main-eq} 
has infinitely many solutions for
$(k,\ell)=(3,2)$ (satisfying $\gcd(n,d)=1$). 
Note that if we
permit negative values of $n$, we must modify this conjecture somewhat to allow
for solutions corresponding to the identities 
$$
\prod_{j=-2 m}^{2 m-1} (2j+1) = \left( \prod_{j=0}^{2m-1} (2j+1) \right)^2
$$
and
$$
\prod_{j=-2m^2-2m}^{2m^2+2m} (2j+1) = \left( (2m+1) \prod_{j=0}^{2m^2+2m-1} (2j+1) \right)^2
$$
where $m$ is a positive integer. 

The literature on equation (\ref{main-eq}) is extensive, dating back to work of
Euler who proved that there are no nontrivial solutions with $(k,\ell) =
(4,2)$.  
It is worth observing that, via an
argument of Granville (unpublished, but reproduced in Laishram and Shorey
\cite{LaSh2}), Erd\H{o}s' conjecture is a consequence of the $abc$-conjecture
of Masser and Oesterl\'e.
Currently, Erd\H{o}s' conjecture has been verified unconditionally only subject to a variety of additional assumptions. By way of example, we now know it to be true if
$d$ is fixed (Marszalek \cite{Mars}), if both $\ell$ and $\omega(d)$ (the number of distinct prime divisors of $d$) are fixed 
(Shorey and Tijdeman \cite{ShTi}), if $P(d)$ (the greatest prime divisor of
$d$) is fixed and $\ell \geq 3$ (Shorey \cite{Sho0}), or if $n$ is fixed and
$\ell \geq 7$ (Shorey \cite{Sho00}).  In subsequent work, a number of these
results have been refined and, in a number of cases, made completely explicit
(particularly for small values of $k$); 
the interested reader is directed to the fine survey of Shorey \cite{Sho1} for
further details on the literature on this problem.

The papers we have mentioned so far rely upon either elementary arguments in
the spirit of Erd\H{o}s and Selfridge, or upon lower bounds for linear forms in
logarithms (sometimes in conjunction with Diophantine inequalities  resulting
from Pad\'e approximation to binomial functions). 
More recently, we find a number of results that appeal to the modularity of
Galois representations associated to certain Frey-Hellegouarch curves to show
that equation (\ref{main-eq})
has at most finitely many solutions, again under certain
additional constraints. The possibility of this approach is implicit in the
work of Darmon and Granville \cite{DaGr} (where, in Corollary 2.1, the
finiteness of the number of nontrivial solutions to (\ref{main-eq}) is proved
provided $k$ and $\ell$ are both fixed). Explicitly,  via such methods, we find
a complete solution 
of equation (\ref{main-eq}) in case $k=3$ (Gy\H{o}ry \cite{Gy}), $k \in \{ 4, 5
\}$ (Gy\H{o}ry, Hajdu and Saradha \cite{GHP}), $6 \leq k \leq 11$ 
(Bennett, Bruin,
Gy\H{o}ry and Hajdu \cite{BBGH}) and $12 \leq k \leq 34$ (Gy\H{o}ry, Hajdu and
Pint\'er \cite{GHP}). In \cite{BBGH}, it is further proved that (\ref{main-eq})
has at most finitely many nontrivial solutions for all $k \leq 82$.

In this paper, we prove a somewhat weakened version of the 
Erd\H{o}s conjecture, which deals also with negative solutions:
 
\begin{thm}\label{thm:main}
There is an effectively computable absolute constant $k_0$ such that if $k \geq
k_0$ is a positive integer, then any solution in integers to 
equation (\ref{main-eq}) with prime exponent $\ell$  
satisfies  either $y=0$ or $d=0$ or $\ell \le \exp(10^k)$. 
\end{thm}
It follows from Faltings' Theorem that \eqref{main-eq} has 
finitely many solutions with $k \ge k_0$ and $yd \ne 0$.

Our proof of Theorem \ref{thm:main} follows
 very different lines from prior work on this problem, and
we emphasize that it 
bears little resemblance to an earlier result of the authors \cite{BeSi}, where
an analogous finiteness statement for rational points on curves corresponding
to equation (\ref{start}) is deduced. While our starting point shares much in
common with \cite{BBGH}, \cite{BeSi} and \cite{GHP}, in that one is led to
study certain ternary equations with corresponding Frey-Hellegouarch curves,
the information we derive from these equations is 
quite distinct from that previously considered. 
In particular, our proof of Theorem \ref{thm:main} 
makes essential use of a wide array of tools
from arithmetic geometry, analytic number theory and additive combinatorics,
including:
\begin{itemize}
\item The modularity of elliptic curves over $\Q$ due to Wiles,
Breuil, Conrad, Diamond and Taylor.
\item Ribet's level lowering theorem.
\item Known cases of Serre's uniformity conjecture,
due to Mazur, to Bilu, Parent \& Rebolledo,
to Darmon \& Merel, and to Lemos.
\item A version of the large sieve inequality due to Selberg.
\item The prime number theorem for Dirichlet $\mathrm{L}$-functions.
\item Gap principles for exceptional zeros of $\mathrm{L}$-functions
due Siegel and Landau.
\item An explicit version of Roth's theorem on $3$-term arithmetic
progressions.
\item Theorems on short character sums due to Burgess and to 
Graham \& Ringrose.
\end{itemize}


The outline of this paper is as follows. In Section \ref{Residual}, we state
some now standard results deriving from the modularity of elliptic curves. In
Section \ref{sec:Frey}, we detail the correspondence  between solutions to
(\ref{main-eq}), related ternary Diophantine equations, and Frey-Hellegoaurch
elliptic curves. We further discuss why the techniques of
 \cite{BBGH} and \cite{GHP} (which lead to analogues of Theorem \ref{thm:main}
for small values of $k$) will likely fail for all sufficiently large $k$.
Sections \ref{First} and \ref{Closer} contain, respectively, an argument that
guarantees that primes in $(k/2,k]$ necessarily divide $d$ (for a 
solution to (\ref{main-eq}) with $y \ne 0$ and large exponent $\ell$), 
and the 
 consequence of this, that the primes $p \equiv 3 \mod{4}$ in this interval are in fact supersingular for a certain parametrized family of elliptic curves. In Section \ref{Character}, we use 
 this information to construct a (short) character sum that is unusually large, corresponding to each Frey-Hellegouarch curve. Section \ref{sec:PNT} contains an argument, based upon the Prime
 Number Theorem for Dirichlet characters, that ensures the desired conclusion, provided we have suitably many elliptic curves corresponding to our Frey-Hellegouarch curves with extremely smooth conductors.
In Section \ref{Consequences}, we attain a like conclusion, via upper bounds for short character sums and the large sieve, under the assumption that we have a somewhat larger number of rather less smooth conductors. Finally, in Sections \ref{enough} and \ref{Pooh}, we complete the proof of Theorem \ref{thm:main}, by using a variety of sieving arguments to show that our Frey-Hellegouarch curves correspond to 
sufficiently many Dirichlet characters to guarantee that we can appeal to at least one of the results from the preceding sections.

\medskip

We are grateful to Adam Harper, Roger Heath-Brown, Lillian Pierce
 and Trevor Wooley for useful 
conversations.

\section{Residual Representations attached to Elliptic Curves} \label{Residual}

Let $E$ be an elliptic curve defined over $\Q$, with minimal discriminant $\Delta$ and conductor $M$.
For a rational prime $\ell \ge 3$, we denote by 
$$
\overline{\rho}_{E,\ell} \; : \; G_{\mathbb{\Q}} \rightarrow \Aut(E[\ell]) \cong \GL_2(\F_\ell)
$$
the representation describing the action of $G_{\mathbb{Q}}:=\Gal(\overline{\Q}/\Q)$ on the $\ell$-torsion subgroup
$E[\ell]$. Define
\begin{equation}\label{eqn:notSerre}
M_0=M\, \Big/ \!
    \prod_{\substack{q \, \| \, M, \; \text{$q$ prime}\\ \ell \, \mid \, \ord_q(\Delta)}} q,
\end{equation}
where we write $\ord_q (x)$ for the largest power of a prime $q$ dividing a nonzero integer $x$.

The following theorem is a standard consequence of Ribet's level lowering theorem \cite{Ribet-1990}
(stated, for example, in \cite[page 157]{Siksek}). It was originally conditional on the modularity
of elliptic curves over $\Q$, a result that was subsequently proved by Wiles, Breuil, Conrad, Diamond and Taylor (see \cite{Wi} and \cite{BreuilConradDiamondTaylor01}).
Additionally, it is, in fact, a special case of Serre's Modularity Conjecture \cite{SerreDuke}, now a theorem of Khare and Wintenberger (\cite{KW1} and \cite{KW2}).  
\begin{thm}\label{thm:Serre}
If $E[\ell]$ is irreducible then there is a cuspidal newform $f=\sum_{n \ge 1} c_n q^n$
of weight $2$ and level $M_0$ such that $\overline{\rho}_{E,\ell} \sim \overline{\rho}_{f,\lambda}$
where $\lambda \mid \ell $ is a prime of the totally real field $K=\Q(c_1,c_2,\dotsc)$.
\end{thm}
Here, by $\overline{\rho}_{E,\ell} \sim \overline{\rho}_{f,\lambda}$ we mean that, for almost all primes $p$, we have that
$$
a_p (E) \equiv c_p \mod{\lambda}.
$$
In fact, by comparing the traces of Frobenius for $\overline{\rho}_{E,\ell}$
and  $\overline{\rho}_{f,\lambda}$, we can be rather more precise.
\begin{lem}\label{lem:traces}
With notation as in Theorem~\ref{thm:Serre}, let $p$ be a rational prime.
\begin{enumerate}
\item[(i)] if $p \nmid \ell M M_0$ then $a_p(E) \equiv c_p \mod{\lambda}$;
\item[(ii)] if $p \nmid \ell M_0$ and $p \, \| \, M$ then $p+1 \equiv \pm c_p \mod{\lambda}$.
\end{enumerate}
\end{lem}

The following lemma will be invaluable to us.
\begin{lem}\label{lem:ellboundp}
With notation as above, suppose $p \ne \ell$ is a prime with $p\, \| \,M$ and, additionally, $\ell \mid \ord_p(\Delta)$. Then 
$$
\ell \le  (\sqrt{p}+1)^{(M_0+1)/6}.
$$
\end{lem}
\begin{proof}
From \eqref{eqn:notSerre}, we see that $p \nmid M_0$. Thus by Lemma~\ref{lem:traces}
we have
$$
\lambda \mid (p+1 \mp c_p)
$$
and so
$$
\ell \mid \Norm_{K/\Q}(p+1\mp c_p).
$$
As $c_p$ is bounded by $2 \sqrt{p}$ in all the real embeddings of $K$, we have
$$
\ell \; \le \; (p+1+2\sqrt{p})^{[K:\Q]} \; =\; (\sqrt{p}+1)^{2[K:\Q]}.
$$
If we denote the dimension of $S_2^{\mathrm{new}}(M_0)$ by $g^+_0(M_0)$, then 
$[K:\Q] \le g_0^+(M_0)$. By Theorem 2 of Martin \cite{Ma}, we have
\begin{equation}\label{eqn:Martin}
g_0^+(M_0) \le \frac{M_0+1}{12},
\end{equation}
completing the proof.
\end{proof}

It is well-known that if the residual characteristic $\ell$ is 
sufficiently large compared to the level $M_0$ then
$f$ has rational eigenvalues and so corresponds to an
elliptic curve over $F/\Q$. We shall have use of a quantitative
version of this statement due to Kraus \cite{Kraus}.
For a positive integer $n$ let
\begin{equation}\label{eqn:mun}
\mu(n)=n  \prod_{\substack{q \mid n\\ {\text{$q$ prime}}}} \left( 1+ \frac{1}{q} \right).
\end{equation}
Define
$$
F(n)= \left(\sqrt{\frac{\mu(n)}{6}} +1\right)^{2 g_0^+(n)}, \; \; \; G(n)=\left(\sqrt{\frac{\mu(\lcm(n,4))}{6}}+1 \right)^2
$$
and set
$$
H(n)=\max(F(n),G(n)).
$$
The following is Th\'{e}or\`{e}me 4 of \cite{Kraus}.
\begin{thm}[Kraus]\label{thm:Kraus}
With notation as in Theorem~\ref{thm:Serre}, suppose $E$ has full $2$-torsion
and that
\[
\ell> H(M_0).
\]
Then there is an elliptic curve $F/\Q$ having 
full $2$-torsion of conductor $M_0$
such that $\overline{\rho}_{E,\ell}\sim \overline{\rho}_{F,\ell}$.
\end{thm}

\section{Frey-Hellegouarch Curves Associated to \eqref{main-eq}}  \label{sec:Frey}

We shall call a solution $(n,d,k,y,\ell)$ 
of \eqref{main-eq} \textbf{trivial} if
$yd=0$. We shall henceforth restrict our attention to nontrivial
solutions.
In this section, we will show how a nontrivial solution to 
equation (\ref{main-eq}) is simultaneously a solution to many generalized
Fermat equations,
both of signature $(\ell,\ell,\ell)$ and of signature $(\ell,\ell,2)$ (in fact, we can actually derive ternary equations of signature $(\ell,\ell,q)$ 
for values of $q > 2$, but these will not be of interest to us).
The following elementary lemma is an immediate consequence of the coprimality
assumption for equation \eqref{main-eq}.
\begin{lem}\label{lem:gcd}
Let $(n,d,k,y,\ell)$ be a nontrivial
solution to \eqref{main-eq} with $\ell$ prime.
\begin{itemize}
\item[(i)] For $0 \le i <j \le k-1$,
\[
\gcd(n+id,n+jd) \mid (j-i).
\]
\item[(ii)] Let $0 \le i \le k-1$ and let $q\ge k$ be prime. Then
\[
\ell \mid \ord_q(n+id).
\]
\end{itemize}
\end{lem}

Thus we may write
\begin{equation}\label{eqn:Ai}
n+id = A_i \, y_i^{\ell},  \qquad 0 \le i \le k-1,
\end{equation}
where $A_i$ are positive integers divisible only by primes $<k$, whereas $y_i$
are divisible only by primes $\ge k$.

\subsection{Fermat Equations of Signature $(\ell,\ell,\ell)$}
In general, given any integers 
$$
0 \leq i_1 < i_2 < i_3 \leq k-1,
$$
the identity
$$
(i_3-i_2)(n+i_1d)+ (i_1-i_3)(n+i_2d) + (i_2-i_1)(n+i_3d)=0
$$
leads to a ternary Diophantine equations of signature $(\ell,\ell,\ell)$. This provides us with roughly $k^3/6$
generalized Fermat equations to consider. For our purposes, it will be convenient to restrict our attention to indices $(i_1,i_2,i_3)$ in arithmetic progression (of which there are approximately $k^2/4$).
Let 
$$
\cA=\{\; 
(i,j,2j-i) \\; \; : \; \;  i, j, 2j-i\in \{ 0,1,\dotsc,k-1\}, \; \;  i<j\; 
\} 
$$
denote the set of nontrivial $3$-term arithmetic progressions
in the set $\{0,1,\dotsc,k-1\}$.
Associated to any such tuple $\ga=(i,j,2j-i) \in \cA$
is the identity
$$
(n+id)-2(n+jd)+(n+(2j-i)d)=0,
$$
from which we see that $(r,s,t)=(y_i,y_j,y_{2j-i})$ is a solution to the
following generalized Fermat equation of signature $(\ell,\ell,\ell)$:
$$
A_i r^\ell-2 A_j s^\ell+A_{2j-i} t^\ell=0 \, .
$$
We may attach to this solution a Frey-Hellegouarch curve as in Kraus  \cite{Kraus}.
For convenience we let
\begin{equation}\label{eqn:ggcd}
g=\gcd \left( n+id,\; 2(n+jd),\; n+(2j-i)d \right),
\end{equation}
\begin{equation}\label{eqn:abc}
a_\ga=\frac{n+id}{g}, \; \;  b_\ga=\frac{-2(n+jd)}{g} \; \; \mbox{ and } \; \; 
c_\ga=\frac{n+(2j-i)d}{g},
\end{equation}
Our corresponding Frey--Hellegouarch is
$$
E_\ga \; : \; Y^2=X(X-a_\ga)(X+ c_\ga).
$$
\begin{lem}\label{lem:Frey1}
The model $E_\ga$ is minimal and semistable at all odd primes.
Its discriminant is
\[
\Delta_\ga=64 (a_\ga b_\ga c_\ga)^2=\frac{2^8}{g^6} (n+id)^2 (n+jd)^2 (n+(2j-i)d)^2 .
\]
In particular, for any prime $p \ge k$, we have $\ell \mid \ord_p(\Delta_\ga)$.
\end{lem}
\begin{proof}
The first part is a straightforward computation. The second follows
from Lemma~\ref{lem:gcd}.
\end{proof}

\begin{lem}\label{lem:level1}
Let $\ell \ge 7$. Then $\overline{\rho}_{E_\ga,\ell} \sim \overline{\rho}_{f,\lambda}$ where $f$ is a newform of weight $2$ and level $M_\ga$, with
\begin{equation}\label{eqn:levelA_i}
M_\ga \mid 2^8 \cdot A_i A_j A_{2j-i},
\end{equation}
and
\[
M_\ga \le 2^7 \cdot \exp(1.000081 \cdot k).
\]
\end{lem}
\begin{proof}
As $E_\ga$ has full $2$-torsion and $\ell \ge 7$, we know from
the work of Mazur \cite{Mazur} that $E_\ga[\ell]$ is irreducible.
It follows from Theorem~\ref{thm:Serre} that $\overline{\rho}_{E_\ga,\ell} \sim \overline{\rho}_{f,\lambda}$ where $f$ is a newform of weight $2$
and level $M_0$ given by \eqref{eqn:notSerre}.
We write $M_\ga:=M_0$. Equation \eqref{eqn:notSerre} and 
Lemma~\ref{lem:Frey1} ensure that $M_\ga$ satisfies
\eqref{eqn:levelA_i}. Moreover, as the odd part of $M_\ga$
is squarefree, $M_\ga$ divides
$$
2^7 \prod_{\substack{q \le k\\ \text{$q$ prime}}} q \, . 
$$
From Schoenfeld \cite[page 160]{Sc}, we have
\begin{equation} \label{gumby}
\sum_{q \le k} \log{q} \; < \; 1.000081 \cdot k.
\end{equation}
The lemma follows.
\end{proof}

\subsection{Fermat Equations of Signature $(\ell,\ell,2)$} \label{subsec:ll2}
Let
$$
\cI=
\{ \;
(j_1,i_1,i_2,j_2) \; \; : \; \; 
i_1+i_2=j_1+j_2, \; \;  0 \leq  j_1<i_1 \leq i_2<j_2 \leq k-1
\;
\}.
$$
To any fixed quadruple $\ii=(j_1,i_1,i_2,j_2) \in \cI$, we can associate the identity
$$
(n+j_1 d)(n+j_2d)-(n+i_1d)(n+i_2d)=(j_1 j_2 - i_1 i_2) d^2. 
$$
It follows that $(r,s,t)=(y_{j_1} y_{j_2},y_{i_1} y_{i_2},d)$
is a solution to the following generalized Fermat equation
with signature $(\ell,\ell,2)$:
\begin{equation} \label{froggie}
A_{j_1} A_{j_2} \cdot r^\ell-A_{i_1} A_{i_2} \cdot s^\ell=(j_1 j_2 - i_1 i_2)
\cdot t^2.
\end{equation}
Following Bennett and Skinner \cite{BennettSkinner}, solutions to this equation also correspond to Frey-Hellegouarch elliptic
curves defined over $\Q$. To simplify notation, write
\begin{equation}\label{eqn:ABkappa}
A=(n+j_1 d)(n+j_2d), \; \;  B=(n+i_1d)(n+i_2d) \; \; \mbox{ and } \; \;  \kappa=j_1 j_2 - i_1 i_2,
\end{equation}
so that
\begin{equation}\label{eqn:ABkappa2}
A-B=\kappa d^2.
\end{equation}
Let
\[
\cE_\ii \; : \; Y^2=X(X^2+2\kappa d X + \kappa A).
\]
\begin{lem}\label{lem:Frey2}
The model $\cE_\ii$ is minimal and semistable at all primes $p \ge k$ that also satisfy $p \nmid \kappa$.
It has discriminant
$$
\Delta_\ii=-64 \kappa^3 A^2 B.
$$
In particular, for any prime $p \ge k$ with $p \nmid \kappa$, we have $\ell \mid \ord_p(\Delta_\ii)$.
\end{lem}
\begin{proof}
This again follows from a straightforward computation with the help of
Lemma~\ref{lem:gcd}.
\end{proof}

\begin{lem}\label{lem:level2}
Let $\ell \ge 11$. Then $\overline{\rho}_{\cE_\ii,\ell} \sim \overline{\rho}_{f,\lambda}$ where $f$ is a newform of weight $2$ and level $M_\ii$ satisfying
$$
M_\ii \le 2^7 \cdot 3^5 \cdot k^4 \cdot \exp(2.000162 \cdot k).
$$
\end{lem}
\begin{proof}
As $\cE_\ii$ has a rational point of order $2$ and $\ell \ge 11$, we know from
the work of Mazur \cite{Mazur} that $\cE_\ii[\ell]$ is irreducible.
It follows from Theorem~\ref{thm:Serre} that $\overline{\rho}_{\cE_\ii,\ell} \sim \overline{\rho}_{f,\lambda}$ where $f$ is a newform of weight $2$
and level $M_0$ given by \eqref{eqn:notSerre}.
We write $M_\ii:=M_0$. Equation \eqref{eqn:notSerre}, together
with Lemma~\ref{lem:Frey2}, 
ensures that $M_\ii$ divides
$$
2^7 \cdot 3^5 \cdot \kappa^2 \cdot \prod_{\substack{q \le k\\ \text{$q$ prime}}} q^2 \, . 
$$
As $\lvert \kappa \rvert < k^2$, the lemma follows from inequality (\ref{gumby}).
\end{proof}

At this point, it is worth mentioning why the techniques of \cite{BBGH} and \cite{GHP} are apparently insufficient to prove Theorem \ref{thm:main} (yet do allow one to show that equation (\ref{main-eq}) has at most finitely many nontrivial solutions for small values of $k$). Intrinsically, they rely upon the fact that for suitably small $k$, and each possible tuple
$$
{\bf{A}} = \left( \mbox{Rad} (A_0), \mbox{Rad} (A_1),  \ldots, \mbox{Rad} (A_{k-1}) \right)
$$
(here, the $A_i$ are as in (\ref{eqn:Ai}); the number of such tuples depends only upon $k$ and not $\ell$ or $d$), we can find $\ii=(j_1,i_1,i_2,j_2) \in \cI$ such that the corresponding polynomial-exponential equation
\begin{equation} \label{s-unit}
x+y = z^2,
\end{equation}
where $z \in \Q$ and $x, y$ are $S$-units, for 
$$
S = \left\{ p  \mbox{ prime}\; : \; p \mid A_{j_1} A_{j_2} A_{i_1} A_{i_2} (j_1 j_2 - i_1 i_2) \right\},
$$
has only ``trivial'' solutions. As a first step, one applies an argument to guarantee that 
$$
p \mid A_1 A_2 \cdots A_{k-1} \implies p < \tau k,
$$
for certain $\tau \in (0,1]$. That we may take $\tau=1$ is immediate from the definition of $A_i$, while, for example, Lemma \ref{lem:pdividesd} of the next section implies a like result with $\tau = 1/2$. It is not especially difficult to improve this to $\tau=1/3$, but it appears to be quite hard to reduce this significantly. From a result of Erd\H{o}s, Stewart and Tijdeman (see e.g. Theorem 4 of \cite{EST}), the number of solutions to equation (\ref{s-unit}) with $x$ and $y$ rational numbers supported on primes of size at most $\tau k$ exceeds $\exp \left( 3 \frac{\sqrt{\tau k}}{\log k} \right)$ for large enough $k$. Since the number of tuples ${\bf{A}}$ to be treated also grows exponentially in $\tau k$, while the cardinality of $\cI$ is 
$$
\sum_{j=2}^{k-1} (k-j) \left[ j/2 \right] = \frac{k^3}{12} - \frac{k^2}{8} - \frac{k}{12} + \frac{\delta}{8}, \; \; \mbox{ where } \; \; \delta = 
\left\{
\begin{array}{ll}
0 & \mbox{ if $k$ is even,} \\
1 & \mbox{ if $k$ is odd,} \\
\end{array}
\right.
$$
our expectation is that for all sufficiently large $k$, there will correspond to each choice of $\ii \in \cI$ a tuple ${\bf A}$ for which the associated equation of the shape (\ref{s-unit}) has nontrivial solutions.

We will  proceed in a very different direction. Rather than attempting to reduce the problem of treating equation (\ref{main-eq}) to that of solving associated ternary equations (which, as we have noted, is likely to be futile for large $k$), we will, in the next two sections,  instead deduce from a nontrivial solution to (\ref{main-eq}) the existence of a large number of elliptic curves that, on some level, mimic the behaviour of elliptic curves with complex multiplication (despite not possessing this property).

\section{A First Result on Primes $k/2<p \le k$} \label{First}

We begin with an easy lemma that ensures that primes in the interval $(k/2,k]$
fail to divide  $A_0 A_1 \cdots A_{k-1}$ for suitably large $\ell$. This
apparently innocuous result (a version of which first appeared in the proof of
Theorem 1.5 of \cite{BBGH}) is actually the key first step in proving Theorem
\ref{thm:main}.  \begin{lem}\label{lem:pdividesd}
Let $k \ge 10^8$ and suppose that $(n,d,k,y,\ell)$ is a nontrivial solution  to \eqref{main-eq} 
with prime exponent $\ell>\exp(10^k)$. 
Let $p$ be a prime in the range $k/2<p \le k$. Then $p \mid d$.
\end{lem}
\begin{proof}
Suppose that $p \nmid d$. Then $p$ divides at least one and at most two of the 
terms $n+d,n+2d,\dotsc,n+kd$. Suppose first that $p$ divides precisely
one such term, say $ p \mid n+id$. It follows from \eqref{main-eq} that
$$
\ell \mid \ord_p(n+id).
$$
Let $\ga$ be any triple of indices in $\cA$
containing $i$. It follows from Lemma~\ref{lem:Frey1} that $E_\ga$
is semistable at $p$ with multiplicative reduction, and that $\ell \mid \ord_p(\Delta_\ga)$.
Applying Lemma~\ref{lem:ellboundp},  we see that 
$$
\ell \le (\sqrt{p}+1)^{(M_\ga+1)/6} .
$$
Now the bound in Lemma~\ref{lem:level1} for $M_\ga$ contradicts the assumption $\ell> \exp(10^k)$.

If instead $p$ divides divides precisely two terms, say $p \mid n+id$ and $p \mid n+(i+p)d$, then we choose $\ii=(i,i+1,i+p-1,i+p) \in \cI$. Let $A$, $B$, $\kappa$ and $d$ be as in
(\ref{subsec:ll2}). From \eqref{main-eq} and \eqref{eqn:ABkappa}, we have
$$
p \mid A, \; \;  \ell \mid \ord_p(A) \;  \mbox{ and }  \;  p \nmid B.
$$
Equation \eqref{eqn:ABkappa2} thus implies that $p \nmid \kappa$ and so
the model $\cE_\ii$ has multiplicative reduction at $p$.
Applying Lemma~\ref{lem:ellboundp}, we see that
$$
\ell \le (\sqrt{p}+1)^{(M_\ii+1)/6}.
$$
Now the bound in Lemma~\ref{lem:level2} for $M_\ii$ contradicts the assumption $\ell> \exp(10^k)$, completing the proof of Lemma \ref{lem:pdividesd}.
\end{proof}

\section{A Closer Look at the Frey-Hellegouarch Curve $E_\ga$} \label{Closer}

The Frey-Hellegouarch curves $\cE_\ii$ associated to $\ii \in \cI$
have been valuable in proving Lemma~\ref{lem:pdividesd}.
We shall not, however, have further use for them and will  instead focus, here and henceforth, 
solely on the Frey-Hellegouarch curves $E_\ga$ associated to 
the $3$-term arithmetic progressions $\ga \in \cA$.

\begin{lem}\label{lem:Fga}
Let $k \ge 10^8$ and suppose that $(n,d,k,y,\ell)$ is a nontrivial solution  to \eqref{main-eq} 
with $\ell>\exp(10^k)$ prime. 
Let $\ga \in \cA$. Then there is an elliptic curve $F_\ga/\Q$
having full rational $2$-torsion and conductor $M_\ga$ such that $\overline{\rho}_{E_\ga,\ell} \sim \overline{\rho}_{F_\ga,\ell}$.
\end{lem}
\begin{proof}
By Theorem~\ref{thm:Kraus}, it is sufficient to show that $\ell> H(M_\ga)$.
From Tenenbaum \cite{Tenenbaum} (Theorem 9 and the remark following it), we have
$$
\prod_{\substack{q \le k\\ \text{$q$ prime}}} 
\left( 1+ \frac{1}{q} \right) \; \le \;
\exp\left(0.27+\frac{5}{\log{k}}\right) \cdot \log{k}.
$$
As $k \ge 10^8$, we obtain 
$$
\prod_{q \le k} \left( 1+ \frac{1}{q} \right) \le 2 \log{k}.
$$
This together with Lemma~\ref{lem:level1} and its proof, shows that 
$\mu(M_\ga)$ and $\mu(\lcm(M_\ga,4))$ are both bounded by
$$
2^8 \log{k} \cdot \exp(1.000081 \cdot k).
$$
Using the previously cited estimate \eqref{eqn:Martin} to bound $g_0^+(M_\ga)$, we easily deduce that
$H(M_\ga)<\exp(10^k)< \ell$ as required.
\end{proof}
Throughout the remainder of the paper, we maintain the assumption $\ell>\exp(10^k)$.
Further, $F_\ga$ will always denote the elliptic curve associated to $\ga$ by Lemma~\ref{lem:Fga}.

\begin{lem}\label{lem:supersingular}
With notation and assumptions as in Lemma~\ref{lem:Fga},
let $p$ be a prime satisfying $k/2<p \le k$.
Then $p$ is a prime of good reduction for both $E_\ga$ and $F_\ga$,
and we have $a_p(E_\ga)=a_p(F_\ga)$. If, moreover, $p \equiv 3 \mod{4}$,
then 
$a_p(F_\ga)=0$
and hence
$p$ is a prime of supersingular reduction for $F_\ga$.
\end{lem}
\begin{proof}
By Lemma~\ref{lem:pdividesd}, we know that every prime $k/2<p \le k$
divides $d$. As $\gcd(n,d)=1$ we see that $p \nmid (n+id)$ for all $i$.
It follows from Lemma~\ref{lem:Frey1} that $p$ is a prime of good reduction
for $E_\ga$. Since the conductor $M_\ga$ of $F_\ga$ is a divisor 
of the conductor of $E_\ga$ (see equation \eqref{eqn:notSerre}),
it follows that $p$  is a prime of good reduction for both elliptic curves.
Hence, by Lemma~\ref{lem:traces}, we know that $a_p(E_\ga) \equiv a_p(F_\ga) \mod{\ell}$.
By the Hasse--Weil bounds $\lvert a_p(E_\ga) - a_p(F_\ga) \rvert \le 4 \sqrt{k}$, whereby the inequality $\ell>\exp(10^k)$
immediately implies that $a_p(E_\ga)=a_p(F_\ga)$.

Let $g$ be as in \eqref{eqn:ggcd}, so that the reduction of $E_\ga$ modulo $p$ is
$$
\tilde{E}_\ga \; : \; Y^2=X(X-n/g)(X+n/g) \, .
$$
If $p \equiv 3 \mod{4}$, then, as is well-known (see e.g. page 41 of \cite{Koblitz}),
$a_p(E_\ga)=0$ whereby also $a_p(F_\ga)=0$.
\end{proof}

Before we proceed, it is worth remarking that Lemma \ref{lem:supersingular} implies that the elliptic curve $F_\ga$ shares supersingular primes with elliptic curves with 
complex multiplication and $j$-invariant $1728$, in the interval 
$k/2<p \le k$. As we shall later observe, $F_\ga$ cannot itself have complex multiplication. This alone, however, is not enough to imply a contradiction; indeed the curve with model
\begin{equation} \label{example}
E \; \; : \; \; Y^2 = X^3 - X + \prod_{p \leq k} p
\end{equation}
has precisely these properties. On the other hand, if we can deduce the existence of an $\ga \in \cA$ for which the conductor of $F_\ga$ is suitably ``small'' (notice that $E$ in (\ref{example}) has conductor that is exponentially large in $k$), then we can apply an effective version of the Chebotarev density theorem to derive a contradiction for large $k$, solely from $F_\ga$ having a surplus of supersingular primes
in the interval $(k/2,k]$ (see Serre \cite{SerreIHES} and Elkies \cite{Elk} for upper bounds on the number of supersingular primes in intervals, for elliptic curves without complex multiplication, both conditional on the 
Generalized Riemann Hypothesis (GRH) and otherwise). 
As we shall observe in Section \ref{enough}, we can guarantee the existence of an $\ga$ for which the conductor of $F_\ga$ is bounded above by $k^\lambda$ for some absolute positive constant $\lambda$. This is sufficient to contradict the Chebotarev density theorem under GRH, but not unconditionally. If we had an $\ga \in \cA$ for which $F_\ga$ has conductor bounded by $(\log k)^\lambda$, say, then we would have an alternative proof of Theorem \ref{thm:main} via this approach. At present, we are unable to prove the existence of such an $\ga$.

\section{On a Character Sum Associated to $F_\ga$} \label{Character}

Henceforth, $F_\ga$ will denote the elliptic curve over $\Q$ having full 
$2$-torsion and conductor $M_\ga$ attached, via Lemma~\ref{lem:Fga}, to a
$3$-term arithmetic progression $\ga \in \cA$, where $\cA$ corresponds to a nontrivial solution  of \eqref{main-eq}.
For a positive integer $N$, we write 
$N^\mathrm{odd}=N \cdot 2^{-\ord_2 (N)}$ for the odd part of $N$.
As usual, we denote by $\Lambda$ the von Mangoldt function
$$
\Lambda (n) = 
\left\{ 
\begin{array}{ll}
\log p & \mbox{ if  $n=p^k$ for some prime $p$ and integer $k \geq 1$,} \\
0 & \mbox{ otherwise}. \\
\end{array}
\right.
$$
\begin{prop}\label{prop:charsum}
Let $k \ge 2 \times 10^{10}$ and let $\ell>\exp(10^k)$ be prime.
Let $(n,d,k,y,\ell)$ be a nontrivial solution  to equation \eqref{main-eq} and 
suppose that $\ga \in \cA$.
Then there exists a quadratic character $\chi_\ga$ 
that is 
primitive of conductor $N_\ga$ 
such that
\begin{equation}\label{eqn:charsumbasic}
\left|  \sum_{k/2 < m \le k} \chi_\ga(m) \cdot \Lambda(m) \right| \; > \; 0.1239 \, k \; .
\end{equation}
Moreover, we have that $N_\ga^{\mathrm{odd}} \mid M_\ga$ and $N_\ga^{\mathrm{odd}} \ne 1$. 
\end{prop}

\noindent \textbf{Remark.}
After proving Proposition~\ref{prop:charsum},
the key to the proof of Theorem~\ref{thm:main} will be to show,
for $k$ suitably large, that
if $N_\ga^{\mathrm{odd}} \ne 1$ for all $\ga$,  
then there is some $\ga$ for which 
the left-hand side of the inequality \eqref{eqn:charsumbasic}
is much smaller than $0.1239 k$.

\subsection*{Legendre Elliptic Curves}
Let
$\lambda \in \Q \setminus \{0,1\}$ and write
\begin{equation} \label{eqn:Legendre}
F_\lambda \; \; : \; \; Y^2 = X (X-1)(X-\lambda),
\end{equation}
often called a {\it Legendre elliptic curve with parameter $\lambda$}.
For $\ga \in \cA$, the elliptic curve $F_\ga$ has full $2$-torsion,
and hence is a quadratic twist of a Legendre elliptic curve $F_\lambda$,
where there are in fact six possible choices for $\lambda$. Define
$$
\frak{S} = \left\{ -t^2 \; : \; t \in \Q \right\} \cup \left\{ 2 t^2 \; : \; t \in \Q \right\}.
$$
We partition $\cA$ into two disjoint subsets, $\cA^{(I)}$ and $\cA^{(II)}$.
\begin{enumerate}[leftmargin=2.5cm]
\item[$\cA^{(I)}$:] This consists of $\ga \in \cA$ such that
at least one of the $\lambda$-invariants of $F_\ga$  
lies outside $\frak{S}$.
\item[$\cA^{(II)}$:] This consists of $\ga \in \cA$ such that  
every $\lambda$-invariant of $F_\ga$ is in $\frak{S}$.

\end{enumerate}
 The precise construction of the character $\chi_\ga$ in the proof of Proposition~\ref{prop:charsum} depends on whether $\ga$ belongs to $\cA^{(I)}$
or $\cA^{(II)}$, but in either case it is closely related to the
$\lambda$-invariants of $F_\ga$.

We require some preliminary results.
\begin{lem}\label{lem:Legendre}
Let $F/\Q$ be an elliptic curve of conductor $M$, semistable away from $2$ (i.e with $M^\mathrm{odd}$ squarefree), having full rational $2$-torsion.
 Let  $\lambda \in \Q$ be any of the six $\lambda$-invariants of $F$.
Then the following hold.
\begin{enumerate}
\item[(i)] $\ord_p(\lambda)=\ord_p(1-\lambda)=0$ for all odd primes $p$
of good reduction for $F$.
\item[(ii)] Let $\omega \in \{\pm 1, \pm 2\}$ and let $\chi$ be the unique primitive quadratic character of conductor $N$ which satisfies
\begin{equation}\label{eqn:chidef}
\chi(p)=\left( \frac{\omega \cdot \lambda}{p}\right)
\end{equation}
for odd primes $p$ with 
$\ord_p(\lambda)=0$. 
Then $N^\mathrm{odd} \mid M$.
\end{enumerate}
\end{lem}
\begin{proof}
As $F$ has full rational $2$-torsion and is
semistable away from $2$, it has a model of the form
$$
F \; : \; Y^2=X(X-a)(X-b)
$$
where $a$, $b$, $a-b$ are non-zero integers with no odd prime common factors. The primes dividing $M^\mathrm{odd}$ are precisely
the odd primes dividing $ab(a-b)$. Since the six  associated $\lambda$-invariants are
$$
b/a, \; \;  a/b, \; \; (a-b)/a, \; \; a/(a-b), \; \; b/(b-a) \; \mbox{ and } \; (b-a)/b,
$$
the lemma follows immediately.
\end{proof}
\begin{lem}\label{lem:mod4}
Let $p \equiv 3 \mod{4}$ be prime and suppose that $F/\F_p$ is an elliptic curve
of the form
$$
F \; : \; Y^2=X(X-1)(X-\eta^2)
$$
for some $\eta \in \F_p \setminus\{0,1,-1\}$. Then 
$F(\F_p)$ contains a subgroup isomorphic to $\Z/2\Z \times \Z/4\Z$.
\end{lem}
\begin{proof}
Since $F$ has full rational $2$-torsion it is enough to show that $F/\F_p$
has a point of order $4$ or, in other words, that one of the 
three points of order $2$ is $2$-divisible. We know $(a,b) \in F(\F_p)$
is $2$-divisible if $a$, $a-1$ and $a-\eta^2$ are all squares.
Suppose $(1,0)$ is not $2$-divisible. Then $1-\eta^2$ is not
a square. As $p \equiv 3 \mod{4}$ it follows that 
$\eta^2-1$ is a square. Thus the point $(\eta^2,0)$ is $2$-divisible. 
\end{proof}

We are now ready to apply this to the elliptic curves $F_\ga$
that arise from solutions to \eqref{main-eq}.
\begin{lem}\label{lem:kro}
Let $k \ge 10^8$ and suppose that $\ell>\exp(10^k)$ is prime.
Assume that $(n,d,k,y,\ell)$ is a nontrivial solution to equation \eqref{main-eq}.
Let $\ga \in \cA$,
and let $\lambda$ be any of the six $\lambda$-invariants of $F_\ga$.
If $p\equiv 3 \mod{8}$ is a prime in the interval
$k/2<p \le k$, then
$$
\left(\frac{\lambda}{p} \right)=-1.
$$
\end{lem}
\begin{proof}
From Lemma~\ref{lem:supersingular}, we know that 
$p$ is a prime of good supersingular reduction for $F_\ga$.
Lemma~\ref{lem:Legendre} tells us that $\ord_p(\lambda)=\ord_p(1-\lambda)=0$,
whence $p$ is a prime of good reduction for $F_\lambda$.
Now $F_\lambda$ is a quadratic twist of $F_\ga$ and so
must also have supersingular reduction at $p$.
In particular $a_p(F_\lambda)=0$, so that 
$$
\#F_\lambda(\F_p)=p+1 \equiv 4 \mod{8}. 
$$
On the other hand, if we suppose that $\lambda$ is a square modulo $p$, then we know from
Lemma~\ref{lem:mod4} that $8 \mid \# F_\lambda(\F_p)$. The resulting  contradiction completes the proof.
\end{proof}

\subsection*{Proof of Proposition~\ref{prop:charsum} for $\ga \in \cA^{(I)}$}

We are ready to prove Proposition~\ref{prop:charsum} for $\ga \in \cA^{(I)}$.
Fix a $\lambda$-invariant
of $F_\ga$ with $\lambda \not\in \frak{S}$.
Suppose first
that $\lambda=t^2$ or $\lambda=-2t^2$ for some non-zero rational $t$.
By the results of \cite{RamareRumely}, the assumption that
$k \ge 2 \times 10^{10}$ forces the existence of (many) primes
$p \equiv 3 \mod{8}$ in the interval
$k/2<p \le k$. For each such prime, we have $\left(\frac{\lambda}{p}\right)=1$,
contradicting Lemma~\ref{lem:kro}. We may therefore
suppose 
\begin{equation} \label{flup}
\lambda \not\in \{ \pm t^2 \; : \; t \in \Q \} \cup \{ \pm 2t^2 \; : \; t \in \Q \}.
\end{equation}
If $a$ and $m$ are relatively prime integers, we write
$$
\vartheta(X;a,m)=\sum_{\substack{{p \le X}\\{p \equiv a \bmod{m}}}} \log{p} 
$$
for the first Chebychev function associated to the arithmetic progression
$a \bmod{m}$. Here, the sum is over primes $p$.
By \cite{RamareRumely}, using the inequality $k \ge 2 \times 10^{10}$, we have
$$
\sum_{\substack{{k/2< p \le k}\\{p \equiv 3 \bmod{8}}}} \log{p} \;
= \; \vartheta(k;3,8) - \vartheta(k/2;3,8) \; \ge \; (1-3 \varepsilon) \cdot \frac{k}{8} 
$$
where $\varepsilon=0.002811$. From Lemma~\ref{lem:kro}, we thus have
\begin{equation}\label{eqn:presum}
 \sum_{\substack{{k/2< p \le k}\\{p \equiv 3 \bmod{8}}}}
- \left(\frac{\lambda}{p} \right)
 \log{p} \;
 \ge \; (1-3 \varepsilon) \cdot \frac{k}{8} 
\end{equation}

Let $\mu_i$ be the primitive quadratic Dirichlet
characters which on odd primes $p$ away from the support of $\lambda$ are given
by
$$
\mu_1(p)=\left(\frac{\lambda}{p}\right), \; \; 
\mu_2(p)=\left(\frac{-\lambda}{p}\right), \; \; 
\mu_3(p)=\left(\frac{2 \lambda}{p}\right)  \;  \mbox{ and } \;  
\mu_4(p)=\left(\frac{-2\lambda}{p} \right),
$$
and observe that
\[
\mu_1(p)-\mu_2(p)-\mu_3(p)+\mu_4(p)=\begin{cases} 
4\left(\frac{\lambda}{p}\right)  & \text{if $p \equiv 3 \mod{8}$}\\
0 & \text{otherwise}.
\end{cases}
\]
We may thus rewrite inequality \eqref{eqn:presum} as
$$
 \sum_{k/2< p \le k}
\left(-\mu_1(p)+\mu_2(p)+\mu_3(p)-\mu_4(p)\right)
 \log{p} \;
 \ge \; (1-3 \varepsilon) \cdot \frac{k}{2},
$$
whereby there necessarily exists some $i \in \{ 1, 2, 3, 4 \}$ such that 
\begin{equation}\label{eqn:muip}
\left| \sum_{k/2<p\le k} \mu_i (p) \log(p)   \right| \; \ge \; (1-3 \epsilon) \cdot
\frac{k}{8}.
\end{equation}
We let $\chi_\ga=\mu_i$
and write $N_\ga$ for its conductor. From (\ref{flup}), we have $N_\ga^{\mathrm{odd}} \ne 1$.
Moreover, by Lemma~\ref{lem:Legendre} we have $N_\ga^{\mathrm{odd}} \mid M_\ga$.
Finally, the left-hand side of \eqref{eqn:charsumbasic}
agrees with the left-hand side of  \eqref{eqn:muip},  
except on $m=q^r$ where $q$ is prime and $r \ge 2$. 
Thus the difference between the two sums is bounded by
$$
\lvert \psi(k)-\vartheta(k)-\psi(k/2)+\vartheta(k/2) \rvert,
$$
where $\vartheta$ and
$\psi$ are the first and second Chebychev functions. 
From (5.3*) and (5.4*) of Theorem 6* of Schoenfeld \cite{Sc}, we have (\ref{eqn:charsumbasic})
as desired. This completes the proof of Proposition~\ref{prop:charsum}
in Case (I).

\subsection*{Legendre Elliptic Curves Revisited}
Let $\lambda \in \Q \setminus \{0,1\}$,  $F_\lambda$
be as in \eqref{eqn:Legendre} and suppose that $p$ is an odd prime satisfying $\ord_p(\lambda)=\ord_p(1-\lambda)=0$.
We will need to use the $2$-descent homomorphism: 
$$
\Theta_\lambda \; : \; F_\lambda(\F_p) \rightarrow \F_p^*/
{\F_p^*}^2 \times 
\F_p^*/{\F_p^*}^2 \times 
\F_p^*/{\F_p^*}^2, \qquad
\Theta_\lambda(Q)=(\theta_1(Q),\theta_2(Q),\theta_3(Q)).
$$
The kernel of $\Theta_\lambda$ is precisely $2 F_\lambda(\F_p)$.
If $Q\ne (0,0)$ then $\theta_1(Q)=x(Q) {\F_p^*}^2$. If $Q \ne (1,0)$
then $\theta_2(Q)=(x(Q)-1){\F_p^*}^2$. If $Q \ne (\lambda,0)$
then $\theta_3(Q)=(x(Q)-\lambda){\F_p^*}^2$. Moreover
$\theta_1(Q)\theta_2(Q)\theta_3(Q)=1 {\F_p^*}^2$ for all 
$Q \in F_\lambda(\F_p)$, which allows us to compute $\Theta_\lambda$
even for the points of order $2$.

\begin{lem}\label{lem:Fminus1}
Let $F_{-1}$ be as in \eqref{eqn:Legendre}
and $p \equiv 5 \mod{8}$ be prime.
Then $2^3 \, \| \, \# F_{-1}(\F_p)$.
\end{lem}
\begin{proof}
We use the fact that $2$ represents the class of
non-squares in $\F_p^*/{\F_p^*}^2$. The images of the
points of order $2$ under $\Theta_{-1}$ are
\[
\Theta_{-1}(0,0)=(1,1,1),
\qquad
\Theta_{-1}(1,0)=(1,2,2),
\qquad
\Theta_{-1}(-1,0)=(1,2,2).
\]
It follows that only $(0,0)$ is $2$-divisible. We find that
$2(i,1-i)=(0,0)$ (where $i^2=-1$ in $\F_p$). The points of order
$4$ are $(i,1-i)$, $(i,1-i)+(0,0)$, $(i,1-i)+(1,0)$, $(i,1-i)+(-1,0)$.
The images of all of these under $\Theta_{-1}$ have $i {\F_p^*}^2$ as
first coordinate. This is not a square in $\F_p$ (as $p \equiv 5 \mod{8}$) and
hence none of the points of order $4$ are $2$-divisible. It follows
that $2^3 \, \| \, \# F_{-1}(\F_p)$.
\end{proof}

\subsection*{Some Preliminary Results for $\ga \in \cA^{(II)}$}
Let $\ga \in \cA^{(II)}$. The proof of Proposition~\ref{prop:charsum}
in this case
is a little harder and requires some further preparation.
By the definition of $\cA^{(II)}$, every $\lambda$-invariant of $F_\ga$ belongs to the set
$\frak{S}$.
Note that, if $\lambda$ is any of the $\lambda$-invariants of $F_\ga$ and we write
$$
\lambda_1=\lambda, \; \; \lambda_2=1-\lambda
\; \mbox{ and } \;  \lambda_3=(\lambda-1)/\lambda,
$$
then the six $\lambda$-invariants of $F_\ga$ are
precisely $\lambda_i^{\pm 1}$ with $i=1, 2$ and $3$. If we have that 
$\lambda = -t^2$ for some rational number $t$, it follows that necessarily 
there exists a rational number $v$ such that $\lambda_2 = 2 v^2$ (whence
$\lambda_3 = 2 (v/t)^2$). Similarly, if we have $\lambda = 2 t^2$ for $t \in
\Q$, then either $\lambda_2$ or $\lambda_3$ is of the shape $2 v^2$ for
rational $v$. In all cases, renaming if necessary, we deduce the existence of (positive)
rational numbers $t$ and $v$ such that 
\begin{equation} \label{football}
\lambda = 2t^2 \; \; \mbox{ and } \; 1- \lambda = 2 v^2,
\end{equation}
whereby $2t^2+2v^2=1$. 

\begin{lem} 
\label{bad-1}
Let $k \ge 10^8$ and suppose that $\ell>\exp(10^k)$ is prime.
Let $(n,d,k,y,\ell)$ be a nontrivial solution to equation \eqref{main-eq} with corresponding $\cA$.
Let $\ga \in \cA$,
and suppose that $\lambda$, 
one of the six $\lambda$-invariants of $F_\ga$, satisfies \eqref{football}
for positive rational numbers $t$ and $v$.
If $p \equiv 5 \mod{8}$ is prime with $k/2<p \le k$, then $\ord_p(t)=\ord_p(v)=0$ and
$$
\left( \frac{tv}{p} \right)=1.
$$
\end{lem}
\begin{proof}
Fix a prime $p \equiv 5 \mod{8}$ with $k/2<p \le k$. By 
Lemma~\ref{lem:supersingular}, $p$
is a prime of good reduction for both $E_\ga$ and $F_\ga$,
and we have $a_p(E_\ga)=a_p(F_\ga)$. By Lemma~\ref{lem:Legendre},
$\ord_p(\lambda)=\ord_p(1-\lambda)=0$ and so, 
from (\ref{football}), $\ord_p(t)=\ord_p(v)=0$.
From the proof of Lemma~\ref{lem:supersingular},
the reduction of $E_\ga$ modulo $p$ is a quadratic twist
of $F_{-1}$, whereby $a_p(F_\ga)=a_p(E_\ga)=\pm a_p(F_{-1})$. 
On the other hand, $F_\lambda$ is a quadratic twist of $F_\ga$ and so
$a_p(F_\lambda)=\pm a_p(F_{-1})$.
If we consider also the quadratic twist of $F_\lambda$ by $2$
$$
F_\lambda^\prime \; : \; Y^2=X(X-2)(X-2\lambda),
$$
since $2$ is a non-square modulo $p$, 
it follows that $a_p(F_\lambda^\prime)=-a_p(F_\lambda)$.
Thus either $a_p(F_\lambda)=a_p(F_{-1})$ or $a_p(F_\lambda^\prime)=a_p(F_{-1})$.
Since Lemma~\ref{lem:Fminus1} implies that  $2^3 \, \| \, \# F_{-1}(\F_p)$, we may conclude that either
$2^3 \, \| \, \# F_\lambda(\F_p)$ 
or $2^3 \, \| \, \# F_\lambda^\prime(\F_p)$.

Now let $\Theta$ be the $2$-descent map for 
$F_\lambda/\F_p$ as given previously. From \eqref{football},
we find that 
$$
\Theta(0,0)=(2,1,2),
\; \; 
\Theta(1,0)=(1,2,2)
\; \; \mbox{ and } \; \; 
\Theta(\lambda,0)=(2,2,1).
$$
It follows that none of the points of order $2$ are $2$-divisible,
and so $2^3 \nmid \#F_\lambda(\F_p)$. 
Hence $2^3 \, \| \, \# F_\lambda^\prime(\F_p)$.

We denote the $2$-descent map for $F_\lambda^\prime$ by $\Theta^\prime$.
The images of the points of order $2$ in $F_\lambda^\prime$ are
$$
\Theta^\prime(0,0)=(2,2,1),
\; \; 
\Theta^\prime(2,0)=(2,2,1)
\; \; \mbox{ and } \; \; 
\Theta^\prime(2\lambda,0)=(1,1,1).
$$
It follows that only $(2\lambda,0)$ is $2$-divisible. Let $i$ be
any square-root of $-1$ in $\F_p$ and set 
$$
P=\left(4ivt +2\lambda, (128iv^5 - 64iv^3)t - 128v^6 + 96v^4 - 16v^2\right) \; \in \; E(\F_p).
$$
Then $2P=(2\lambda,0)$ and so $P$ is a point of order $4$. 
Writing $\Theta^\prime=(\theta_1^\prime,\theta_2^\prime,\theta_3^\prime)$, we have that $\theta_3^\prime(P)=4itv \cdot {\F_p^*}^2$. 
Suppose
$$
\left(\frac{4itv}{p} \right)=1.
$$
Then $\theta_3^\prime(P)=1$ and so $\Theta^\prime(P)=(1,1,1)$ or $(2,2,1)$
(recall that the product of the entries is a square). Hence
either $\Theta^\prime(P)=(1,1,1)$ or $\Theta^\prime(P+(0,0))=(1,1,1)$.
It follows that one of the points of order $4$ is $2$-divisible
and so $F_\lambda^\prime(\F_p)$ contains a subgroup isomorphic to $\Z/2\Z \times \Z/8\Z$,
contradicting the fact that $2^3 \, \| \, \# F_\lambda^\prime(\F_p)$. 
We therefore have that
$$
\left(\frac{4itv}{p} \right)=-1
$$
and hence the fact that $i$ is a non-square modulo $p$ completes the proof.
\end{proof}

\subsection*{Proof of Proposition~\ref{prop:charsum} for $\ga \in \cA^{(II)}$}
By an easy modification of our earlier argument, but now using Lemma~\ref{bad-1} in place of Lemma~\ref{lem:kro}, the inequality \eqref{eqn:charsumbasic}
is satisfied, where now $\chi_\ga$ is
a primitive quadratic character which for odd primes 
away from the support of $tv$ is given by
\[
\chi_\ga(p)=\left(\frac{\omega \cdot tv}{p}\right)
\]
for some $\omega \in \{\pm 1,\pm 2\}$ that depends only on $\ga$.
Again we write $N_\ga$ for the conductor
of $\chi_\ga$.

We would like to show that $N_\ga^{\mathrm{odd}} \mid M_\ga$.
We may choose a model for $F_\ga$ of the form
$Y^2=X(X-a)(X-b)$
where $a$, $b$ and $a-b$ are non-zero integers, with no odd prime common factors, and we have
$$
2t^2=\lambda=b/a \; \; \mbox{ and } \; \;  2v^2=1-\lambda=(a-b)/a.
$$
Thus
the odd primes appearing in the support of $\omega \cdot tv$
are primes dividing $a$, $b$ or $a-b$. As $\chi_\ga$
is quadratic, $N_\ga^{\mathrm{odd}}$, the odd part of its conductor, 
is squarefree.
On the other hand, the primes dividing $M_\ga^{\mathrm{odd}}$ are precisely
the odd primes dividing $ab(a-b)$, whereby $N_\ga^{\mathrm{odd}} \mid M_\ga$
as required.

Finally, we must prove that $N_\ga^{\mathrm{odd}} \ne 1$.
Suppose $N_\ga^{\mathrm{odd}}=1$. Then $tv=\pm \alpha^2$
or $tv=\pm 2 \alpha^2$ for some positive rational $\alpha$.
We have chosen $t$ and $v$ positive, whereby necessarily $tv=\alpha^2$ or $tv=2\alpha^2$.
Write $t=T/U$ and $v=V/U$ where, without loss of generality,
$T$, $V$ and $U$ are positive integers with $\gcd(U,V,T)=1$. Then, from (\ref{football}), 
$$
2T^2+2V^2=U^2
$$
and hence $T$ and $V$ are odd and coprime, while $U \equiv 2 \mod{4}$.
In particular $2 \mid \ord_2(tv)$ and so we may conclude that $tv=\alpha^2$.
It follows that $TV$ is a positive integer square and hence, since $T$ and $V$ are coprime
and positive, each is itself an integer square, say $T=T_0^2$ and $V=V_0^2$,
where $T_0$ and $V_0$ are positive.
Writing $U=2 U_0$, we thus have
$$
T_0^4+V_0^4=2 U_0^2,
$$
whereby, from a classical descent argument,  $T_0=V_0=U_0=1$, 
and so $\lambda = 1/2$. In particular $F_\ga$ is isomorphic (possibly
over a quadratic extension) to the elliptic
curve $Y^2=X(X-1)(X-1/2)$ with $j$-invariant $1728$ and
complex multiplication by $\Z[i]$. It follows that $F_\ga$ has complex
multiplication and
hence the image of $\overline{\rho}_{F_\ga,\ell}$ is contained
in the normalizer of a Cartan subgroup of $\GL_2(\F_\ell)$. As
$\overline{\rho}_{E_\ga,\ell} \sim \overline{\rho}_{F_\ga,\ell}$
the same is trivially true for $\overline{\rho}_{E_\ga,\ell}$.
It follows
from the work of Lemos \cite{Lemos} (building on 
the results of Darmon and Merel \cite{DM} and of Bilu, Parent and
Rebolledo \cite{BPR}) that $E_\ga$ also has complex
multiplication.  
If we let $a=a_\ga$, $b=b_\ga$ and $c=c_\ga$ be as in \eqref{eqn:abc},
we find that the
$j$-invariant of $E_\ga$ is 
$$
j=2^8 \frac{(a^2-bc)^3}{a^2 b^2 c^2}.
$$
Since $E_\ga$ has complex
multiplication, $j$ is integral. The fact that $a$, $b$ and $c$ are coprime thus implies that each of $a$, $b$ and $c$
is not divisible by odd primes. As $a+b+c=0$, we quickly deduce that
two out of $a$, $b$ and $c$ are equal. 
If $a=b$ or $c=b$ then
$$
n+id=-2(n+jd) \; \;  \mbox{ or } \; \;  n+(2j-i)d=-2(n+jd)
$$
which imply that
$$
3n=-(2j+i)d \; \;  \mbox{ or } \; \; 3n=(i-4j)d.
$$
Since $\gcd(n,d)=1$, it follows that $d \mid 3$, contradicting
Lemma~\ref{lem:pdividesd}.
We thus have $a=c$ and so $n+id=n+(2j-i)d$, whence $d=0$.
The resulting contradiction completes the proof of Proposition~\ref{prop:charsum}.

\section{The Prime Number Theorem} \label{sec:PNT}

Henceforth we fix a nontrivial solution $(n,d,k,y,\ell)$ to equation \eqref{main-eq} (with corresponding $\cA$), and suppose that
$\ell$ and $k$ satisfy the assumptions of Proposition~\ref{prop:charsum}. By this proposition, 
$N_\ga^{\mathrm{odd}} \ne 1$,
and therefore $\chi_\ga$ is \textbf{nontrivial} for each $\ga \in \cA$,
a fact that will be crucial in obtaining a bound for $k$.

\medskip

We shall make use of the Prime Number Theorem for Dirichlet 
characters. Let us begin by defining what we mean by {\it exceptional conductors} and {\it exceptional zeros} for Dirichlet $L$-functions; here we combine Theorems 5.26 and 5.28 of \cite{IK}.

\begin{prop}
There exists an effectively computable absolute constant $c^*>0$ such that 
the following hold.
\begin{enumerate}
\item[(i)] If $\chi_1$ and $\chi_2$ are distinct real, primitive quadratic
characters of conductor $N_1$ and $N_2$, respectively, with associated
$L$-functions $L(s,\chi_1)$ and $L(s,\chi_2)$ having real zeros
$\beta_{\chi_1}$ and $\beta_{\chi_2}$, respectively, then
\begin{equation} \label{two}
\min \{ \beta_{\chi_1}, \beta_{\chi_2} \} < 1 - \frac{3c^*}{\log (N_1 N_2)}.
\end{equation}
\item[(ii)]  If $\chi$ is any primitive, quadratic character of conductor $N$, then $L(s,\chi)$ has at most a single real zero $\beta_\chi$ with
\begin{equation} \label{one}
1 - \frac{c^*}{\log N} < \beta_\chi < 1.
\end{equation}
If such a zero exists, then $\chi$ is necessarily real and $\beta_\chi$ is a simple zero. We term $\beta_\chi$ an \textbf{exceptional zero} and $N$ an 
\textbf{exceptional conductor}.
\end{enumerate}
\end{prop}

From this, if $N_1 < N_2$ are two exceptional conductors, with corresponding
exceptional zeros $\beta_{\chi_1}$ and $\beta_{\chi_2}$,  then, combining
(\ref{two}) and (\ref{one}), 
$$
1 - \frac{c^*}{\log N_1} < \min \{ \beta_{\chi_1}, \beta_{\chi_2} \}  < 1 - \frac{3c^*}{\log (N_1 N_2)},
$$
and so 
\begin{equation} \label{growth}
N_2 > N_1^2.
\end{equation}

The following quite explicit version of the Prime Number Theorem for Dirichlet characters  is Theorem 5.27 of \cite{IK}.
\begin{thm} \label{thm:PNT}
Let $\chi$ be a primitive Dirichlet character of conductor $N$. Then
\begin{equation} \label{PNT}
\sum_{m \le X} \chi(m) \Lambda(m)=\delta_\chi X-\frac{X^{\beta_\chi}}{\beta_\chi}
+O\left(X \exp\left(\frac{-c \log{X}}{\sqrt{\log{X}}+\log{N}}\right) \cdot
(\log{N})^4\right) .
\end{equation}
Here $\delta_\chi=0$ unless $\chi$ is trivial in which case $\delta_\chi=1$.
Moreover, $c>0$ is an absolute effective constant, and the implied
constant is absolute. Also $\beta_\chi$ denotes the exceptional zero if 
present, otherwise the term $-X^{\beta_\chi}/\beta_\chi$ is to be omitted. 
\end{thm}

It is worth observing at this point that the ``error term'' here is actually
smaller than the main term (so that the statement in non-trivial), only for
suitably small conductor $N$, relative to the interval of summation $X$; i.e.
only when $\log N \ll \log^\kappa X$ for some $\kappa < 1$. We wish to apply
this result to characters of conductor roughly $N_\ga$, over an interval of
length $k/2$. The difficulty we encounter is that, {\it a priori}, the $N_\ga$
can be as large as $e^k$ and, even on average, are of size that grows
polynomially in $k$. Further, the potential presence of an exceptional
(Siegel-Landau) zero $\beta_\chi$ additionally complicates matters, even when
we have $N_\ga$ much smaller than $k$, as the term on the right-hand side of
(\ref{PNT}) corresponding to $\beta_\chi$ can, potentially, be very close to
$k$ in size. If, however, we are able to show that we can find sufficiently
many $\ga$ for which $N_\ga$ is ``tiny'', we can use the fact that exceptional
conductors are rare (as quantified in inequality (\ref{growth}), a ``repulsion
principle'' due to Landau), to  reach the desired conclusion :

\begin{prop}\label{lem:tinycond}
Let us suppose that $0 < c_1 < 1$ is fixed and, further, that there is a subset $\cD$ of $\cA$ such that the following
hold :
\begin{enumerate}
\item[(i)] 
$P(N_\ga) \ne 
P(N_\ga^\prime)$ whenever $\ga \ne \ga^\prime$
belong to $\cD$; 
\item[(ii)] $P(N_\ga)<(\log{k})^{1-c_1}$ for all $\ga \in \cD$;
\item[(iii)] 
\begin{equation}\label{eqn:Precipbound}
\sum_{\ga \in \cD} \frac{1}{P(N_\ga)} \ge 0.166.
\end{equation}
\end{enumerate}
Then there exists an effectively computable constant $k_1$,  depending only upon $c_1$, such that $k \leq k_1$.
\end{prop}

We will later apply this proposition with $c_1 = 10^{-4}$. The constant $0.166$ is chosen so that, in our argument,  we have enough progressions $\ga$ to guarantee that either one corresponds to a non-exceptional conductor, or, through appeal to (\ref{growth}), that the smallest exceptional conductor $N_\ga$ we encounter satisfies $N_\ga \leq 400000$, contradicting work of Platt \cite{Pl}.

To prove Proposition \ref{lem:tinycond}, it is convenient for us to be able to deduce an explicit upper bound upon $N_\ga$, given one for $P(N_\ga)$.
\begin{lem}\label{lem:condbound}
Let $N$ the conductor of a quadratic character, and let $P(N)$
be the largest prime factor of $N$. Then $P(N)>  0.94 \log{N}$.
\end{lem}
\begin{proof}
We can write $N=2^\kappa N_1$, where $N_1$ is squarefree and $\kappa \in \{ 0, 1, 2 \}$. 
Then
$$
\log{N} \le \kappa \log{2}+\sum_{p \le P(N)} \log{p}  < \kappa \log 2 + 1.000081 P(N),
$$
via  work of Schoenfeld \cite[page 160]{Sc}. We thus have
$$
\frac{P(N)}{\log (N)} > 0.9999 \left( 1 - \frac{\kappa \log 2}{\log (N)} \right).
$$
The desired result is then immediate if $\kappa =0$ (i.e. unless $4 \mid N$). If $\kappa=1$, we have the claimed inequality, unless $N \leq 105932$, while, for $\kappa=2$, the conclusion follows for all $N \geq 1.2 \times 10^{10}$. 
A (relatively) short computation, checking values of $N \equiv 4 \mod{8}$ up to $105932$ and $N \equiv 8 \mod{16}$ to $1.2 \times 10^{10}$ with, in each case, the odd part of $N$ squarefree, completes the proof;  the minimum value of $P(N)/\log (N)$ is attained at $N=24$.
\end{proof}

\begin{proof}[Proof of Proposition  \ref{lem:tinycond}]
Suppose there is some $\ga \in \cD$ such that 
the character $\chi_\ga$ 
is non-exceptional. 
By assumption (ii)
and Lemma~\ref{lem:condbound}, 
 $\log{N_\ga} < 1.07  \, (\log{k})^{1-c_1}$.
Applying  Theorem \ref{thm:PNT}, we have
$$
\sum_{k/2<m \le k} \chi_\ga(m) \Lambda(m)=
O\left(k \exp\left(-c^\prime (\log{k})^{c_1}\right) \cdot
(\log{k})^4\right) \, ,
$$
for some effectively computable positive constant $c^\prime$, contradicting \eqref{eqn:charsumbasic} for $k$ sufficiently large.

We may therefore suppose that $\chi_\ga$ is exceptional
for every $\ga \in \cD$.
We obtain, from assumption (i),
a sequence of exceptional conductors
$$
N_1<N_2<\cdots<N_s
$$
where $s=\#\cD$. From inequality (\ref{growth}),
$N_j> N_1^{2^{j-1}}$, whence, via Lemma~\ref{lem:condbound},
$$
P(N_j)> 0.94 \cdot 2^{j-1} \log{N_1},
$$
for each $j$.
By assumption \eqref{eqn:Precipbound},
$$
0.166 \le \sum_{j=1}^s \frac{1}{P(N_j)}<\frac{2.13}{\log{N_1}},
$$
whereby
$$
N_1 \le 373743,
$$
contradicting work of Platt \cite{Pl}, which rules out exceptional zeros
corresponding to Dirichlet characters, for every conductor smaller than
$400000$.  
\end{proof}

\section{Consequences of having enough characters $\chi_\ga$ with smooth, small conductors} \label{Consequences}

In the previous section, we stated a result (Proposition \ref{lem:tinycond}) that guarantees an effective upper bound upon $k$, provided we have suitably many $\ga$ with $P(N_\ga)$ ``tiny'', i.e. with $N_\ga$ very smooth.
In this section, we will show that, in fact, we can reach the same conclusion if we have a (potentially) much larger number of somewhat less smooth conductors corresponding to $\ga \in \cA$.
\begin{prop}\label{lem:cB}
Suppose that $c_2>10$ is a constant and that there exists a subset $\cB \subset \cA$
such that 
\begin{enumerate}
\item[(i)] $\# \cB > 17 \log{k}$;
\item[(ii)] for every distinct pair $\ga$, $\ga^\prime \in \cB$ 
we have $\chi_\ga \ne \chi_{\ga^\prime}$; 
\item[(iii)]  
$P(N_\ga) \le k^{7/16}$ for all $\ga \in \cB$;
\item[(iv)] $N_\ga<k^{c_2}$.
\end{enumerate}
Then there is an effectively computable constant $k_2$,  depending only upon $c_2$, such that $k \leq k_2$.
\end{prop}

Here, the constants $17$ and $7/16$ can be slightly sharpened, but this is not
of great importance for our argument.

The proof of Proposition \ref{lem:cB} relies upon a combination of ingredients,
including the large sieve and upper bounds for character sums over short
intervals. We begin with the latter. 

\subsection{Character Sums over Short Intervals}

We shall need a standard theorem on short character sums
to a smooth modulus, a variant 
of some results of Graham and Ringrose \cite{GrRi}. Specifically, we will appeal to \cite[Theorem 12.13]{IK}.
\begin{thm}\label{thm:GR1}
Let $\pi_i$ be characters of conductor $q_i$, for $1 \le i \le r$. 
Write $q=q_1$ and suppose  that $q>1$ is
squarefree with $\gcd(q,q_2 q_3 \cdots q_{r})=1$.
Suppose, moreover, 
that $\pi_1 $ is primitive. 
Then, for $R \ge R_0$ where
$$
R_0 = \max(q_2,\dots,q_{r},q^{1/4}) \, q^{5/4},
$$
we have
$$
\left\lvert
\sum_{M<m \le M+R} \pi_1\cdots \pi_{r-1} \pi_r(m) 
\right\rvert
\le
4 R \cdot \left(\tau(q)^{r^2}/q \right)^{2^{-r}},
$$
where $\tau(q)$ is the number of divisors of $q$.
\end{thm}

We will prove the following.
\begin{prop}\label{lem:GR2}
Let $c_2>0$ be a constant. Then there exist effectively computable positive
constants $k_3$ and 
$c_3$, each depending only on $c_2$,
such that the following holds.
Let $k \ge k_3$ be an integer and suppose that
$\chi_1$ and $\chi_2$ are \textbf{distinct}
primitive quadratic characters modulo $N_1$ and $N_2$, respectively,
where the $N_i$ satisfy   
\begin{equation} \label{conditions}
P(N_i) \le k^{7/16} \; \; \mbox{ and } \; \; N_i \leq k^{c_2}, \; \; \mbox{ for } \; i \in \{ 1, 2 \}.
\end{equation}
Then
\begin{equation}\label{eqn:GR1}
\left\lvert
\sum_{k/2<m \le k} \chi_1(m) \chi_2(m) 
\right\rvert
\; \le \;
 k^{1- c_3} \, .
\end{equation}
\end{prop}
\begin{proof}
Let $\chi=\chi_1 \chi_2$ and write $M=\lcm(N_1,N_2)$ for the conductor of $\chi$. We can thus rewrite $\chi=\eta \psi$
where $\eta$ is \textbf{primitive} of conductor $M_1$ and $\psi$ is 
\textbf{principal}
of conductor $M_2$ with $M=M_1 M_2$ and $\gcd(M_1,M_2)=1$. 
As $\eta$ is quadratic, we see that $M_1^{\mathrm{odd}}$ is squarefree.
Clearly, $M_2 \mid \gcd(N_1,N_2)$, and so $M_2^{\mathrm{odd}}$ is also
squarefree.
From \eqref{conditions}, 
\begin{equation}\label{eqn:conditions2}
P(M)\le k^{7/16} \; \; \mbox{ and } \; \;  M \le k^{2c_2}.
\end{equation}
We shall consider two cases, according to whether $M_1 \ge 8 k^{7/32}$
or $M_1< 8 k^{7/32}$.

\vskip 5mm

\noindent \textbf{Case 1}. Suppose first that 
\begin{equation}\label{eqn:M1big}
M_1 \ge 8 k^{7/32},
\end{equation}
so that
$$
M_1^{\mathrm{odd}} \ge k^{7/32}.
$$
We can write
$$
\eta=\pi_1\dotsc \pi_s \; \; \mbox{ and } \; \;  \psi=\pi_{s+1}\dotsc \pi_{r},
$$
where $\pi_i$ is primitive of modulus $q_i$ for $i=1,\dotsc,s$
and principal of modulus $q_i$ for $i={s+1},\dotsc, r$. Moreover,
the $q_i$  (which could be composite) may be chosen to satisfy
\begin{enumerate}
\item[(a)] $q_1 q_2 \dotsc q_s=M_1$ 
and $q_{s+1} q_{s+2} \cdots q_r=M_2$,
\item[(b)] $q_1 \mid M_1^{\mathrm{odd}}$ and 
so $\gcd(q_1,q_2 q_3 \cdots q_{r})=1$,
\item[(c)]  $k^{7/32} \le q_i \le k^{7/16}$, for $i=1, \dotsc, s-1$
and $i=s+1,\dotsc,r-1$,
\item[(d)] $1<q_r \le k^{7/16}$ and
\item[(e)] $s \ge 1$, and if $s>1$ then $1 \le q_s \le k^{7/16}$. 
\end{enumerate}
%
Now, from property (c) and \eqref{eqn:conditions2},
$$
r-2   \; \le \;   
\log{M}/\log(k^{7/32})
\; < \;  10 c_2 \, , 
$$
whence $r < 10c_2+2$.
In the notation of Theorem~\ref{thm:GR1}, we have that
$$
R_0 \le k^{7/16} \cdot (k^{7/16})^{5/4} \le k^{63/64} < k/2.
$$
Notice here that, at least in this argument, we cannot replace the exponent $7/16$ in (\ref{conditions}) with one larger than $4/9$.
 
We will now apply Theorem~\ref{thm:GR1}. 
Let $q=q_1$ and note that we have
(see e.g. page 334 of \cite{IK})
$$
\tau(q) \le q^{1/\log\log{3q}},
$$
for all $q \geq 1$.
As $q \ge k^{7/32}$
and $r <10 c_2+2$, we see that for $k$ suitably large, 
$$
\tau(q)^{r^2} < q^{1/2}.
$$
Appealing to Theorem~\ref{thm:GR1}, we thus have
$$
\left\lvert
\sum_{k/2<m \le k} \chi_1(m) \chi_2(m) 
\right\rvert \le \frac{2k}{q^{1/2^{r+1}}},
$$
whence inequality \eqref{eqn:GR1} follows from $q \ge k^{7/32}$ and $r < 10 c_2+2$. Explicitly, we may take $c_3 = 2^{-10c_2-6}$.
This completes the proof of Proposition \ref{lem:GR2} in Case 1.

\vskip 5mm

\noindent \textbf{Case 2}. Next, suppose instead that 
$$
M_1 < 8 k^{7/32}.
$$
Since $\chi_1$ and $\chi_2$ are distinct, it follows 
that $\chi=\chi_1 \chi_2$ is not principal, and so 
$$
\left|
\sum_{k/2<m \le k} \chi_1(m) \chi_2(m) 
\right| < M=M_1 M_2. 
$$
To complete the proof of \eqref{eqn:GR1}, we may thus certainly suppose that
$$
M_2>k^{3/4}.
$$
Write $\mu$ for the M\"{o}bius function, and recall that
$$
\sum_{d \mid n} \mu(d)=\begin{cases}
1 & \mbox{ if } n=1, \\
0 &  \mbox{ if } n>1.
\end{cases}
$$
Now we can write
\[
\begin{split}
\sum_{k/2<m \le k} \chi_1(m) \chi_2(m) & = 
\sum_{k/2<m \le k} \eta(m) \psi(m) \\
&= 
 \sum_{\substack{{k/2< m \le k}\\{\gcd(m,M_2)=1}}} \eta(m)\\
&= 
 \sum_{k/2< m \le k} \eta(m) \sum_{d \mid \gcd(m,M_2)} \mu(d)\\
&= \sum_{d \mid M_2} \sum_{k/2<nd \le k} \eta(nd) \mu(d) \\
&= \sum_{d \mid M_2} \eta(d) \mu(d)  \sum_{k/(2d) <n \le k/d} \eta(n) \\
\end{split}
\]
As $\eta$ is non-principal and has conductor $M_1<8 k^{7/32}$, we have
$$
\left|
\sum_{k/(2d) <n \le k/d} \eta(n) 
\right| < M_1< 8 k^{7/32} \, .
$$
Thus
\[
\left|
\sum_{k/2<m \le k} \chi_1(m) \chi_2(m) 
\right| < \tau(M_2) \cdot 8 k^{7/32}\le M_2^{1/\log\log{3M_2}} \cdot 8 k^{7/32}. 
\]
The proof is complete for $k$ sufficiently large as
$k^{3/4}<M_2<k^{c_2}$.
\end{proof}

\subsection{Proof of Proposition~\ref{lem:cB}: The Large Sieve}

We make use
of the following inequality of Bombieri (Proposition 1 of \cite{Bombieri}, attributed there to Selberg).
\begin{thm} \label{thm:Bombieri}
If $\xx$, $\yy_1,\dotsc,\yy_m$ are vectors in an inner product
space then
\[
\sum_{i=1}^m \lvert \xx \cdot \yy_i \rvert^2 \le 
\lVert \xx \rVert^2 \cdot \max_{1 \le i \le m} \left\{
\sum_{j=1}^m \lvert \yy_i \cdot \yy_j \rvert
\right\} \, .
\]
\end{thm}

\vskip2ex

In view of \eqref{eqn:charsumbasic},
to prove Proposition \ref{lem:cB},  it clearly suffices to show that
\begin{equation}\label{eqn:hope2}
\frac{1}{\#\cB} \sum_{\ga \in \cB}
\left\lvert \sum_{k/2 < m \le k} \chi_\ga (m) \cdot \Lambda(m)
\right\rvert^{2} \le \varpi \cdot k^2,
\end{equation}
for $k$ sufficiently large, where $\varpi=0.1239^2$. 

Let $\xx=( \Lambda(m) )_{k/2 < m \le k}$ and,
for each $\ga \in \cB$,  choose corresponding 
$\yy_\ga=(\chi_\ga(m))_{k/2<m \le k}$ so that the desired inequality \eqref{eqn:hope2} can be rewritten as
\begin{equation}\label{eqn:hope3}
\frac{1}{\# \cB} \sum_{\ga \in \cB} \lvert \xx \cdot \yy_\ga \rvert^2
\le \varpi \cdot k^2 \, .
\end{equation}
Applying the large sieve (Theorem \ref{thm:Bombieri}), we  have
\begin{equation} \label{flippy}
\frac{1}{\# \cB} \sum_{\ga \in \cB} \lvert \xx \cdot \yy_\ga \rvert^2
\le 
\lVert \xx \rVert^2 \cdot \max_{\ga \in \cB} \left\{
\frac{1}{\#\cB} 
\sum_{\ga^\prime \in \cB} \lvert \yy_\ga \cdot \yy_{\ga^\prime} \rvert 
\right\} \, .
\end{equation}

Let us begin by noting that
\[
\begin{split}
\lVert \xx \rVert^2 & =\sum_{k/2 < m \le k} \Lambda(m)^2\\
& \le \log{k} \sum_{k/2 <m \le k} \Lambda(m)\\
& 
=  \frac{k \log{k}}{2}+ O(k),
\end{split}
\]
from the Prime Number Theorem.
Further, for each $\ga \in \cB$,
we have
$$
\left| \yy_\ga \cdot \yy_\ga \right| 
\leq \frac{k+1}{2}.
$$
As $\#\cB \ge 17 \log{k}$ (assumption (i)), it follows that 
$$
\frac{\left| \yy_\ga \cdot \yy_\ga \right|}{\# \cB} \le 
\frac{k+1}{34 \log{k}}.
$$
Next, we would like to estimate $\yy_\ga \cdot \yy_{\ga^\prime}$
for $\ga \ne \ga^\prime$ belonging to $\cB$.
Assumptions (ii), (iii), (iv) ensure that $\chi_\ga$, $\chi_{\ga^\prime}$ 
satisfy the conditions of Proposition~\ref{lem:GR2}, which gives
\[
\lvert \yy_\ga \cdot \yy_{\ga^\prime} \rvert
=\left| \sum_{k/2<m<k} \chi_1(m) \chi_2(m) \right| \; \le \; k^{1-c_3} \, .
\] 
Hence, from (\ref{flippy}),
\begin{equation} \label{mission}
\begin{split}
\frac{1}{\# \cB} \sum_{\ga \in \cB} \lvert \xx \cdot \yy_\ga \rvert^2
& \le 
\lVert \xx \rVert^2 \cdot
\max_{\ga \in \cB} \left\{
\frac{\lvert \yy_\ga \cdot \yy_\ga \rvert}{\#\cB} +
\max_{\ga^\prime \ne \ga} 
\lvert \yy_\ga \cdot \yy_{\ga^\prime} \rvert 
\right\} \\
& \le \left( \frac{k \log{k}}{2}+ O(k) \right) \cdot
\left(\frac{k+1}{34 \log{k}}+ k^{1-c_3}
 \right)  \\
& =\frac{k^2}{68} \cdot \left(1+o(1) \right).
\end{split}
\end{equation}
As $1/68 <  \varpi^2$, we have inequality \eqref{eqn:hope3}, as desired, for $k$ suitably large.
This completes the proof of Proposition~\ref{lem:cB}.

\section{Generating Enough Characters} \label{enough}

We now wish to sieve the set $\cA$ carefully, hoping to guarantee the existence of suitably many corresponding characters $\chi_\ga$ with conductors smooth enough and small enough to enable us to employ either Proposition \ref{lem:tinycond} or Proposition \ref{lem:cB}.
There are (at least) two approaches we can take here to find a reasonable quantity of smooth characters, both dependent upon leaving a positive proportion of elements in $\cA$ after application of our sieve. We could, for example, appeal to a theorem of Varnavides \cite{Varnavides} which guarantees that a set of positive density in $\{ 0, 1, \ldots, k-1 \}$ contains $\gg k^2$ nontrivial $3$-term arithmetic progressions, and then average over these progressions. Instead, we will rely upon an explicit version of a theorem of Roth on $3$-term arithmetic progressions, together with an old argument of Erd\H{o}s. An apparent (small) advantage of this approach is that it will lead to explicit and reasonably small values for $c_2$ in Proposition \ref{lem:cB}.
We begin by stating 
\begin{thm}[Roth]\label{thm:Roth}
Let $0<\delta<1$. Then there exists a positive constant $K_0 (\delta)$ such if $k \geq K_0 (\delta)$ and  $J \subset \{0,1,\dotsc,k-1 \}$ with
$\#J \ge \delta k$, then there is at least one
nontrivial $3$-term arithmetic progressions in $J$, i.e. there exist integers $0 \leq i < j$ such that $i, j$ and $2j-i$ all belong to $J$.
\end{thm}
Note here that, following work of Rahman \cite{Rah}, for example,  we may take
\begin{equation} \label{rahmen}
K_0 (\delta) = \exp (\exp (132 \log (2) \cdot \delta^{-1} ) ).
\end{equation}

Let us define our index set
$I=\{0,1,\dotsc,k-1 \}$ and 
recall that $\cA$ is the set of $3$-term arithmetic progressions
$(i,j,2j-i)$
in $I$, i.e. the set of integer triples $(i,j,2j-i)$, satisfying $0 \leq i<j$ and $2j-i < k$.
For a prime $p$, write
$$
I_p=\{ i \in I \; : \;  p \mid (n+id) \}, 
$$
so that
$$
\# I_p= \delta_p \left(\frac{k}{p}+\theta_p \right)
$$
where $\lvert \theta_p \rvert<1$ and 
$$
\delta_p=
\begin{cases}
1 & \mbox{ if } p \nmid d \\
0 & \mbox{ if } p \mid d
\end{cases} \, .
$$

We will now use Theorem \ref{thm:Roth}, together with an elementary argument of Erd\H{o}s, to find an element of $\ga \in \cA$ with corresponding conductor $N_\ga$ that is smooth, small, and coprime to a given ``thin'' set of primes. We will do this in completely explicit form to provide an indication of the size of the constants involved here (and in particular to demonstrate an admissible value for $c_2$ in Proposition \ref{lem:cB}).

\begin{prop}\label{lem:S}
Let us suppose that 
\begin{equation} \label{big}
k \geq \exp ( \exp ( 10^6 ))
\end{equation}
is an integer and 
that $S \subset [1,k]$ is a set of primes satisfying
\begin{equation}\label{eqn:sumS}
\sum_{p \in S} \frac{1}{p}< 0.17.
\end{equation}
Then there exists an $\ga \in \cA$ satisfying the following:
\begin{enumerate}
\item[(I)] $p \nmid N_\ga$ for $p \in S$;
\item[(II)] $P(N_\ga) \le k^{7/16}$;
\item[(III)] $N_\ga$ is not divisible by primes in the range
$((\log{k})^{1-10^{-4}},10^4 \log{k}]$;
\item[(IV)] $N_\ga<k^{418}$.
\end{enumerate}
\end{prop}
\begin{proof}
Suppose that $k$ satisfies (\ref{big}).
Let us define  $T$ to be the set of primes in the interval $( k^{7/16},k ]$, $U$ to be the primes in the 
interval $((\log{k})^{1-10^{-4}},10^4 \log{k}]$ and set
$$
J=I\setminus \bigcup_{p \in S \cup T \cup U} I_p.
$$
Notice that if $\ga=(i,j,2j-i)$ is an arithmetic progression in $J$,
then $(n+id)$, $(n+jd)$ and $(n+(2j-i)d)$ are each not divisible
by any prime $p$ in $S$, $T$ or $U$.
By \eqref{eqn:levelA_i} and Proposition~\ref{prop:charsum}, the conductor $N_\ga$ therefore 
satisfies (I), (II) and (III).

Our initial goal will be to show that the set $J$ has positive density in $I$.
Note that
$$
\# \bigcup_{p \in S \cup T \cup U} I_p  \le \sum_{p \in S} \# I_p
+\sum_{p \in T} \# I_p + \sum_{p \in U} \# I_p\, . 
$$
Now
$$
\sum_{p \in T} \# I_p =\sum_{k^{7/16} < p \le k} \delta_{p} \left(\frac{k}{p} +\theta_p\right)
$$
and hence we have
$$
\sum_{p \in T} \# I_p  < k 
\sum_{k^{7/16} < p \le k/2} \frac{1}{p}
+ \frac{ 0.6 k}{\log k},
$$
where we have used the fact that $\delta_p=0$ for all $k/2 < p \leq k$,  Theorem 1 of Rosser and Schoenfeld \cite{RoSc}, which yields the inequalities
$$
\frac{x}{\log x} \left( 1 + \frac{1}{2 \log x} \right) < \pi (x) < \frac{x}{\log x} \left( 1 + \frac{3}{2 \log x} \right),
$$
provided $x \geq 59$, and (\ref{big}).
From Theorem 5 of Rosser and Schoenfeld \cite{RoSc}, we have
$$
\left| \sum_{p \le x} \frac{1}{p} - \log\log{x}- \tau \right| < \frac{1}{2 \log^2 x},
$$
valid for $x \geq 286$, where $\tau$ is an absolute constant (explicitly, $\tau = 0.26149 \ldots$),
and hence
$$
\sum_{k^{7/16} < p \le k/2} \frac{1}{p} <  \log (16/7) + \log \left( 1 - \frac{\log 2}{\log k} \right) + \frac{1}{2 \log^2 (k/2)} + \frac{128}{49 \log^2 k}.
$$
From (\ref{big}), we thus have
$$
\sum_{k^{7/16} < p \le k/2} \frac{1}{p} <  \log (16/7) - \frac{0.6}{\log k}
$$
and hence
$$
\sum_{p \in T} \# I_p \le \log(16/7) \cdot k.
$$
Moreover, 
$$
\sum_{p \in U} \frac{1}{p}  < \log\left(\frac{\log\log{k}+\log{10^4}}{(1-10^{-4}) \log\log{k}} \right)
+\frac{1}{\log^2 \left( (\log{k})^{1-10^{-4}} \right)}
$$
and so, from (\ref{big}),
$$
\sum_{p \in U} \frac{1}{p}  < \log \left(1/(1-10^{-4}) \right) + \frac{5 \log 10}{\log \log k},
$$
whence
$$
\sum_{p \in U} \# I_p \le \log \left(1/(1-10^{-4}) \right) k + \frac{5 \log (10) \, k}{\log \log k} + 10^4 \log k.
$$

From  \eqref{eqn:sumS}, we have, crudely,
$$
\sum_{p \in S} \# I_p \le 0.17 \, k+\pi (k) < 0.17 k + \frac{1.1 k}{\log k}.
$$
Thus
$$
\# \bigcup_{p \in S \cup T \cup U} I_p  \le (\log(16/7)+\log(1/(1-10^{-4}))+0.17)k
+ \frac{12k}{\log\log{k}} 
$$
and hence, from (\ref{big}), we have 
$$
\# \bigcup_{p \in S \cup T \cup U} I_p 
< 0.9968  \, k.
$$
It follows that
$$
\# J=\# I - \# \bigcup_{p \in S \cup T \cup U} I_p > 0.0032 k,
$$
so that, in particular, $J$ is nonempty (and, as noted earlier, possesses the property that any arithmetic progression $\ga=(i,j,2j-i)$ in $J$ has corresponding $N_\ga$ satisfying (I), (II) and (III)). From Theorem \ref{thm:Roth}, it is immediate that there exist nontrivial $3$-term arithmetic progressions $\ga$ in $J$; it remains to show that at least one of them has property (IV), i.e. satisfies $N_\ga \leq k^{418}$.

We now follow a classic argument of Erd\H{o}s (see e.g. Lemma 3 of \cite{Er55}, or, in the context of arithmetic progressions, display equation (3.6) of \cite{LaSh2}), defining a set
$J_1 \subset J$, obtained by deleting from $J$, for each prime $p \leq k$, an index $i_p$ with the property that $\ord_p (A_{i_p})$ is maximal. 
It follows that
$$
\#J_1 > 0.0032 k - \pi (k) > 0.00319 k
$$
and, more importantly for our purposes, that
$$
\prod_{i \in J_1} A_i \mid (k-1)!.
$$
Since no prime $p \geq k^{7/16}$ divides any of these $A_i$,  Stirling's formula (see e.g. \cite{Strom} for a suitably explicit version) thus implies that
$$
\prod_{i \in J_1} A_i \leq \sqrt{2 \pi (k-1)} ((k-1)/e)^{k-1} e^{1/(12(k-1))} \prod_{k^{7/16} < p \leq k} p^{-\ord_p((k-1)!)} \, .
$$
Now
\[
\log \left(\prod_{k^{7/16} < p \leq k} p^{\ord_p((k-1)!)} \right)
\ge
\sum_{{k^{7/16} < p \leq k}} \left(\frac{k-1}{p}-1 \right) \log{p}
\ge
\frac{9}{16} k \log{k}-5k
\]
using Theorem 5 of \cite{Tenenbaum}, Theorem 6 of \cite{RoSc} and our
assumption \eqref{big}.
Hence, after a little work, 
$$
\prod_{i \in J_1} A_i < k^{0.44 k}.
$$
It follows, if we define $J_2 \subset J_1$ to be the set of indices $i \in J_1$ with the property that $A_i \leq k^{139}$, that $\# J_2 > 0.00001 k$. Checking that in (\ref{rahmen}) we have
$$
K_0(10^{-5}) < \exp ( \exp ( 10^6 )), 
$$
we may thus apply Theorem \ref{thm:Roth} (Roth's theorem), to deduce the existence of a nontrivial $3$-term arithmetic progression of indices $\ga=(i,j,2j-i)$ in $J_2$.  By \eqref{eqn:levelA_i} 
$$
N_\ga \leq 2^8 \cdot  A_i A_j A_{2j-i} \leq 2^8 \cdot (k^{139})^3 < k^{418}.
$$
This concludes the proof of Proposition \ref{lem:S}.
\end{proof}

\section{Proof of Theorem~\ref{thm:main}} \label{Pooh}

We are now ready to prove Theorem \ref{thm:main}. 
To begin, note that there exists a non-empty subset $\cB \subset \cA$
satisfying
\begin{enumerate}
\item[(i)] $P(N_\ga) \ne P(N_{\ga^\prime})$
whenever $\ga \ne \ga^\prime$ in $\cB$;
\item[(ii)] $P(N_\ga) \le k^{7/16}$ for all $\ga \in \cB$;
\item[(iii)] $N_\ga$ is not divisible by primes in the range
$[(\log{k})^{1-10^{-4}},10^4 \log{k}]$, for all $\ga \in \cB$;
\item[(iv)] $N_\ga<k^{418}$
for all $\ga \in \cB$.
\end{enumerate}
Indeed to generate such a $\cB$ with one element, we may simply 
apply Proposition~\ref{lem:S} with $S=\emptyset$. Now let
$\cB$ be a \textbf{maximal} nonempty subset of $\cA$ satisfying (i)--(iv).
If $\#\cB>17 \log{k}$,  then $k$ is effectively bounded by Proposition~\ref{lem:cB}.
We may thus suppose that $\# \cB \le 17 \log{k}$.
Assume first that
$$
\sum_{\ga \in \cB} \frac{1}{P(N_\ga)}<0.17.
$$
It follows, if we let $S=\{P(N_\ga) \; : \; \ga \in \cB\}$, that $S$ satisfies
\eqref{eqn:sumS}. Proposition~\ref{lem:S} thus yields the existence of
some $\ga \in \cA$ that satisfies (ii), (iii), (iv) and,
moreover, has the property that $N_\ga$ is not divisible by any prime in $S$.
Thus $P(N_\ga) \ne P(N_\ga^\prime)$ for $\ga^\prime \in \cB$.
Now the set $\cB^\prime=\cB \cup \{ \ga\}$ is strictly larger
than $\cB$ and satisfies conditions (i)--(iv), contradicting the maximality of
$\cB$.  

We may thus assume that
$$
\sum_{\ga \in \cB} \frac{1}{P(N_\ga)} \ge 0.17.
$$
Define
$$
\cC=\{\ga \in \cB \; : \; P(N_\ga) > 10^4 \log{k}\}
$$ 
and
$$
\cD=\{ \ga \in \cB \; : \; P(N_\ga) < (\log{k})^{1-10^{-4}} \} .
$$
Then, by condition (iii), $\cB$ is the disjoint union of $\cC$ and $\cD$. It follows that
$$
\sum_{\ga \in \cC} \frac{1}{P(N_\ga)} \le \frac{\#\cC}{10^4 \log{k}}
\le \frac{\# \cB}{10^4 \log{k}} \le \frac{17 \log{k}}{10^4 \log{k}}=0.0017,
$$
whereby
$$
\sum_{\ga \in \cD} \frac{1}{P(N_\ga)} \ge 0.1683.
$$
We now apply Proposition~\ref{lem:tinycond} with $c_1 = 10^{-4}$ to deduce that $k$ is bounded.
This completes the proof of Theorem~\ref{thm:main}.

\section{Concluding remarks} \label{Closing}

Much of the literature on (\ref{main-eq}) has, in fact, dealt with the somewhat more general equation
\begin{equation} \label{lamb}
n ( n+d) \cdots (n+ (k-1)d ) = b y^{\ell},
\end{equation}
where $b$ is an integer, all of whose prime factors are bounded above by $k$. The arguments we have presented here do not permit us to treat quite such a general 
situation, but can be extended to handle equation (\ref{lamb})  where $P(b)$, the greater prime factor of $b$, is at most $\tau k$, for $\tau < 1/2$. 

While we have given our results in Section \ref{Character} on characters attached to nontrivial solutions to (\ref{main-eq}) only  for large values of $k$, analogous statements are readily obtained for smaller $k$.
These provide us with a way to prove that the number of nontrivial solutions to (\ref{main-eq}) is finite that is much more computationally efficient  than that described in \cite{BBGH}. Since the lower bound upon $k$  in Theorem \ref{thm:main} is so large, however, there is little chance we can treat all the remaining cases $k \leq k_0$ by such an approach, without the introduction of  fundamentally new ideas.


\end{document}